\RequirePackage{amsmath}
\documentclass{kluwer}

\usepackage{amsfonts}

\usepackage{amssymb}

\input{xy}
\xyoption{all}

\newcommand{\Cbar}{{\ensuremath{\overline{C}}}}
\newcommand{\Chat}{{\ensuremath{\hat{C}}}}
\newcommand{\pihat}{\ensuremath{\hat{\pi}}}
\newcommand{\spihat}{\ensuremath{\underline{\hat{\pi}}}}
\newcommand{\Crep}{{\ensuremath{\check{C}}}}
\newcommand{\smod}{\ensuremath{\underline{\textup{mod}}\,}}

\newcommand{\Hom}{\ensuremath{\textup{Hom}\,}}
\newcommand{\sHom}{\ensuremath{\underline{\textup{Hom}}\,}}
\newcommand{\Ext}{\ensuremath{\textup{Ext}}}
\newcommand{\End}{\ensuremath{\textup{End}\,}}
\newcommand{\sEnd}{\ensuremath{\underline{\textup{End}\,}}}
\newcommand{\ind}{\ensuremath{\textup{ind}\,}}
\newcommand{\add}{\ensuremath{\textup{add}\,}}
\newcommand{\iadd}{\ensuremath{\textup{iadd}\,}}

\newcommand{\pd}{\ensuremath{\textup{pd}\,}}
\newcommand{\id}{\ensuremath{\textup{id}\,}}
\newcommand{\isomorphe}{\ensuremath{\cong}}
\newcommand{\CA}{{\ensuremath{\mathcal{C}_A}}}
\newcommand{\CC}{{\ensuremath{\mathcal{C}_C}}}
\newcommand{\DA}{{\ensuremath{\mathcal{D}^b(\textup{mod}\, A)}}}
\newcommand{\DC}{{\ensuremath{\mathcal{D}^b(\textup{mod}\, C)}}}

\newcommand{\HD}{\ensuremath{\Hom_\DA(\oplus_{i\in\mathbb{Z}} F^i T,-)}}
\newcommand{\HDX}[1]{\ensuremath{\Hom_\DA(\oplus_{i\in\mathbb{Z}} F^i T,{#1})}}

\newcommand{\zG}{\ensuremath{\Gamma}}
\newcommand{\zS}{\ensuremath{\Sigma}}
\newcommand{\zO}{\ensuremath{\Omega}}
\newcommand{\za}{\ensuremath{\alpha}}
\newcommand{\zb}{\ensuremath{\beta}}
\newcommand{\ze}{\ensuremath{\epsilon}}
\newcommand{\zd}{\ensuremath{\delta}}
\newcommand{\zg}{\ensuremath{\gamma}}
\newcommand{\zl}{\ensuremath{\lambda}}
\newcommand{\zL}{\ensuremath{\Lambda}}

\newtheorem{thm}{Theorem}[section]
\newtheorem{prop}[thm]{Proposition}

\newtheorem{cor}[thm]{Corollary}
\newtheorem{lem}[thm]{Lemma}
\newtheorem{example}[thm]{Example}

\begin{document}
                  
\begin{article}
\begin{opening} 
  \title{On the Galois coverings of a cluster-tilted algebra
  \thanks{The first two authors are partially supported by the NSERC of
  Canada and the university of Sherbrooke, the second author is
  partially supported by Bishop's university and the third author is
  partially supported by the NSF grant DMS-0700358 and the University
  of Massachusetts at Amherst}} 
  \author{Ibrahim \surname{Assem} }
\runningtitle{On the Galois coverings of a cluster-tilted algebra}
\institute{Universit\'e de Sherbrooke}
  \author{Thomas \surname{Br\"ustle}}
\institute{Universit\'e de Sherbrooke and Bishop's University}
  \author{Ralf \surname{Schiffler}}
\institute{University of Massachusetts at Amherst}
\date{\today}

\begin{abstract}
We study  the module category of a certain Galois covering of a
cluster-tilted algebra which we call the cluster repetitive
algebra. Our main result compares the module categories 
of the cluster repetitive algebra of a tilted algebra $C$ and the
repetitive algebra of $C$, in the sense of Hughes and Waschb\"usch.
\end{abstract}
\end{opening}  

\setcounter{section}{-1}

\begin{section}{Introduction}
The cluster category was introduced in 
\cite{BMRRT} and also in \cite{CCS} for type $\mathbb{A}$, as a
categorical model to understand better the  cluster algebras of Fomin
and Zelevinsky \cite{FZ}. It is a quotient of the bounded derived
category $\DA$ of the finitely generated modules over a finite
dimensional hereditary algebra  $A$.  It was then natural to consider
the endomorphism algebras of tilting objects in the cluster
category. Such algebras are called cluster-tilted, and have been the
subject of several investigations since their introduction in
\cite{BMR1,CCS}, see, for instance
\cite{BMR2,CCS2,KR,ABS,ABS2,BFPT}. In particular, it was shown in
\cite{ABS} that the cluster-tilted algebras are trivial extensions of
tilted algebras by a certain bimodule. 

Now, the class of trivial extensions of tilted algebras by the minimal
injective cogenerator has been extensively investigated. They play an
important r\^ole in the classification results for self-injective
algebras. In this study, one of the essential tools is the repetitive
algebra, introduced by Hughes and Waschb\"usch in \cite{HW}. In
previous works \cite{ABST1,ABST2}, we have related the cluster category
and the $m$-cluster category to the repetitive algebra of a hereditary
algebra. 

Our initial motivation in this paper is different. Given a tilted
algebra $C$, we wish to relate the trivial extension $T(C)$ of $C$ by
its minimal injective cogenerator $DC$ and the corresponding
cluster-tilted algebra $\tilde C$. Doing so has been difficult to
achieve directly, so we decided to work instead with certain  Galois
coverings of these two algebras, the repetitive algebra $\Chat$ of $C$,
which is a covering of $T(C)$, and the algebra $\Crep$ constructed in a
similar manner starting from $\tilde C$, which we call the cluster
repetitive algebra. 

Before stating our main theorem, we recall from \cite{BMR1} that, if $\tilde T$ is a tilting object in the cluster category $\CA$, and $\tilde C=\End_\CA\tilde T$, then the functor $\Hom_\CA(\tilde T, -):\CA\to \textup{mod}\, \tilde C$ induces an equivalence $
\CA/\iadd(\tau\tilde T)\isomorphe \textup{mod}\,\tilde C$, where
$\iadd(\tau\tilde T)$ is the ideal 
 consisting of all morphisms which factor
through a direct sum of summands of the Auslander-Reiten translate
$\tau\tilde T$ of $\tilde T$. Our main theorem says that this functor lifts
to a functor $\textup{mod}\, \Chat \to \textup{mod}\, \Crep$ which
satisfies a similar condition. Namely, we give a different realisation
of the cluster category, using only the tilted algebra $C$, which we
denote as $\CC$, then construct two functors $\phi: \textup{mod}\,
\Chat \to \textup{mod}\, \Crep$ and $\pihat: \textup{mod}\, \Chat \to
\CC$ as well as an ideal $\mathcal{J}$ of $\textup{mod}\, \Chat$
which satisfy the properties stated in the following theorem. 

\begin{thm} Let $C$ be a tilted algebra. Then there is a commutative diagram of dense functors
\[
\xymatrix@C=80pt@R=40pt{
\textup{mod}\, \Chat \ar[r]^\phi \ar[d]_\pihat & \textup{mod}\, \Crep
\ar[d]^{G_\zl} \\
\CC \ar[r]^{\Hom_{\CC}(\pihat C, -)} & \textup{mod}\, \tilde C
}
\]
 where $G_\zl:\textup{mod}\, \Crep\to \textup{mod}\, \tilde C$ is the
 push-down functor associated to the covering $\Crep\to \tilde
 C$. Moreover, $\phi$ is full and induces an equivalence of categories
 $\textup{mod}\, \Chat/\mathcal{J} \isomorphe \textup{mod}\,  \Crep$.
\end{thm}

 Note that the functor $G_\zl$
is always dense: this is not true of the push-down functor
$\textup{mod}\,\Chat\to \textup{mod}\,T(C)$ (see, for instance,
\cite{AS}).

As a consequence of this theorem, we are able to relate the Auslander-Reiten 
quivers of $\Chat$ and $\Crep$, this yields the required relation
between $T(C)$ and $\tilde C$.

The paper is organised as follows. After  brief preliminaries, we start by introducing the notion of cluster repetitive algebra and study its most elementary properties in section \ref{sect 1}. In section \ref{sect 2}, we relate the module category of $\Crep$ to the bounded derived category $\DA$, and show that $\textup{mod}\, \Crep$ is equivalent to a quotient of $\DA$ by a certain ideal. Section \ref{sect 3} is devoted to the proof of our main theorem. Finally, in section \ref{sect 4}, motivated by the need to bring down this information to $\textup{mod}\, \tilde C$, we compute a fundamental domain for $\textup{mod}\, \tilde C$ inside $\textup{mod}\, \Crep$, and show that such a domain lies entirely inside a certain finite dimensional quotient of $\Crep$, which we call the cluster duplicated algebra.

\end{section}


\begin{section}{The cluster repetitive algebra }\label{sect 1}

\begin{subsection}{Notation}\label{sect 1.1}
Throughout this  paper, all algebras  are basic locally finite dimensional algebras over an algebraically closed field $k$. For an algebra $C$, we denote  by $\textup{mod}\, C$ the
  category of finitely generated right $C$-modules and by $\ind C$ a
  full subcategory of $\textup{mod}\, C$ consisting of exactly one
  representative from each isomorphism class of indecomposable
  modules. When we speak about a $C$-module (or an indecomposable
  $C$-module), we always mean implicitly that it belongs to
  $\textup{mod}\, C$ (or to $\ind C$, respectively). Also, all
  subcategories of $\textup{mod}\, C$ are full and so are identified
  with their object classes. Given a subcategory $\mathcal{C}$ of
  $\textup{mod}\, C$, we sometimes write $M\in \mathcal{C}$ to express
  that $M$ is an object in $\mathcal{C}$. We denote by $\add
  \mathcal{C}$ the full subcategory of $\textup{mod}\, C$ with objects
  the finite direct sums of modules in $\mathcal{C}$ and, if $M$ is a
  module, we abbreviate $\add\{M\}$ as $\add M$. 

Following \cite{BoG}, we sometimes consider equivalently an algebra $C$ as a locally bounded $k$-category, in which the object class $C_0$ is a complete set $\{e_i\}_i$ of primitive orthogonal idempotents of $C$, and the group of morphisms from $e_i$ to $e_j$ is $e_iCe_j$. We denote the projective (or the injective) dimension of a module $M$ as $\pd M$ (or $\id M$, respectively).  The global dimension of $C$ is denoted by
  gl.dim.$C$. Finally, we
  denote by $\zG(\textup{mod}\, C)$ the Auslander-Reiten quiver of an
  algebra $C$, and by $\tau_C= D\, Tr$, $\tau^{-1}_C=Tr\, D$ its
  Auslander-Reiten translations. For further definitions and facts
  needed on $\textup{mod}\, C$ or $\zG(\textup{mod}\, C)$, we refer
  the reader to  \cite{ASS,ARS}.

\end{subsection}

\begin{subsection}{Cluster-tilted algebras}\label{sect 1.2}
Let $A$ be a  finite dimensional hereditary algebra. The \emph{cluster category}  $\CA$ of $A$
is defined as follows. Let $F$ be the automorphism of
$\DA$ defined as the composition
$\tau^{-1}_{\DA} [1]$, where $\tau^{-1}_{\DA}$ is the Auslander-Reiten
translation in $\DA$ and $[1]$ is the shift functor. Then $\CA$ is the
orbit category $\DA/F$, that is,  the objects of $\CA$ are the
$F$-orbits $\tilde X=(F^iX)_{i\in \mathbf{Z}}$, where $X\in\DA$, and
the set of morphisms from $\tilde X=(F^iX)_{i\in \mathbf{Z}}$ to
$\tilde Y=(F^iY)_{i\in \mathbf{Z}}$ is 
\[\Hom_{\CA}(\tilde X,\tilde Y) = \bigoplus_{i\in \mathbf{Z}}
\Hom_{\DA}(X,F^iY).
\]
It is shown in \cite{BMRRT,K} that $\CA$ is a triangulated
 category with almost split triangles. Furthermore, the projection  $\pi:\DA\to \CA$ is a
 functor of triangulated categories and commutes with the
 Auslander-Reiten translations, see \cite{BMRRT}.

An object $\tilde T$ in $\CA$ is called a \emph{tilting object}
 provided $\Ext^1_{\CA}(\tilde T,\tilde T)=0$ and the number of
 isomorphism classes of indecomposable summands of $\tilde T$ equals
the rank of the Grothendieck group $\textup{K}_0(A)$. The endomorphism algebra 
$B =\End_{\CA}(\tilde T)$ is then called a \emph{cluster-tilted
 algebra}. The functor $\Hom_{\CA}(\tilde T, -):\CA\to \textup{mod}\, B$
 induces an equivalence 
\[\CA/\iadd(\tau_\CA \tilde T) \isomorphe \textup{mod}\, B,\]
where $\tau_\CA$ is the Auslander-Reiten translation in $\CA$ and
 $\iadd(\tau_\CA \tilde T)$ is the ideal of $\CA$ consisting of all
 morphisms which factor through objects of  $\add(\tau_\CA \tilde T)$. Also, the above equivalence commutes with the
 Auslander-Reiten translations in both categories, see
\cite{BMR1}.

Let $B=\End_\CA\tilde T$ be a cluster-tilted algebra. It is shown in \cite{BMRRT} that we may suppose without loss of generality that the object $\tilde T=(F^iT)_{i\in\mathbb{Z}}$ is such that  $T\in \DA$ is an $A$-module. In this case, the algebra $C=\textup{End}_AT$ is tilted, the trivial extension $\tilde C=C\ltimes \Ext^2_C(DC,C) $ is cluster-tilted and, conversely, any cluster-tilted algebra is of this form, see \cite{ABS}. We also need the following easy lemma.

\begin{lem}\label{lem 1.2} Let $T$ be an $A$-module such that $\tilde T=(F^iT)_{i\in\mathbf{Z}}$ is a tilting object in $\CA$, then
\[\Hom_{\DA}(T,\tau F^i T)=0,
\]
for all $i\in\mathbb{Z}$.
\end{lem}  

\begin{pf} This follows from 
\[ \oplus_{i\in\mathbb{Z}}\Hom_\DA(T,\tau F^i T) \isomorphe \Hom_\CA(\pi T,\tau \pi T)=0,\]
because $\tilde T=\pi T$ is tilting in $\CA$.
\qed
\end{pf}  

\end{subsection}

\begin{subsection}{The cluster repetitive algebra }\label{sect 1.3}
Let $C$ be a tilted algebra. We define the \emph{cluster repetitive algebra} to be the following locally finite dimensional algebra without identity
\[
\Crep  \quad = \quad \left[
  \begin{array}{cccccccccc}
  \ddots &&&\ 0\ \\
&\ C_{-1}\ \\
&E_0&\ C_0\ \\
&&E_{1}&\ C_1\ \\
&\ 0\ &&&\ddots
  \end{array}
 \right] \]
where matrices have only finitely many non-zero entries, $C_i=C$ 
and $E_i =\Ext^2_C(DC,C) $ for all $i\in \mathbb{Z}$, all the remaining
entries are zero and the multiplication is induced from that of $C$,
the $C$-$C$-bimodule structure of $\Ext^2_C(DC,C) $ and the zero map
$\Ext^2_C(DC,C) \otimes_C \Ext^2_C(DC,C) \to 0$. The identity maps
$C_i\to C_{i-1}$, $E_i\to E_{i-1}$ induce an automorphism $\varphi$ of
$\Crep$. The orbit category $\Crep/<\varphi>$ inherits from $\Crep$
the structure of a $k$-algebra and is easily seen to be isomorphic to
the cluster-tilted algebra $\tilde C=C\ltimes \Ext^2_C(DC,C) $. The
projection functor $G:\Crep\to \tilde C$ is thus a Galois covering
with infinite cyclic group generated by $\varphi$, see \cite{G}. We
denote be $G_\zl:\textup{mod}\, \Crep\to \textup{mod}\, \tilde C$ the
associated push-down functor. We need another description of $\Crep$.

\begin{lem}\label{lem 1.3}
Let $T$ be a tilting $A$-module, and $C=\textup{End}_AT.$ Then
\[ \Crep\isomorphe\End_\DA(\oplus_{i\in\mathbb{Z}}F^i T).
\]\ 
\end{lem}

\begin{pf} As a $k$-vector space, we have
\[ \End_\DA(\oplus_{i\in\mathbb{Z}}F^i T)
= \oplus_{i,j\in \mathbb{Z}} \Hom_\DA(F^iT,F^jT).
\]
But $ \Hom_\DA(F^iT,F^jT)=0$ unless $i=j$ or $i=j-1$ since $T\in \textup{mod}\, A$. Moreover, $\Hom_\DA(F^iT,F^iT)=\Hom_A(T,T)=C$ and $\Hom_\DA(F^iT,F^{i+1}T)=\Hom_\DA(T,FT)=\Ext^2_C(DC,C) $, where the last isomorphism follows from \cite[Theorem 3.4]{ABS}.
\qed
\end{pf}  

\end{subsection}

\begin{subsection}{The quiver of $\Crep$}\label{sect 1.4}
The quiver $Q_{\Crep}$ of $\Crep$ is easily deduced from the quiver $Q_C$ of $C$. Let $\{e_1,e_2,\ldots,e_n\}$ be a complete set of primitive orthogonal idempotents of $C$, then $\{e_{\ell,i}\mid 1\le\ell\le n, i\in \mathbb{Z}\}$ is a complete set of primitive orthogonal idempotents of $\Crep$. We have moreover
\[ e_{\ell,i}\ \Crep\  e_{h,j} \isomorphe
\left\{
\begin{array}{cl}
e_\ell\  C \ e_h &\textup{if } i=j\\
e_\ell\ \Ext^2_C(DC,C) \  e_h &\textup{if } i=j+1\\
0& \textup{otherwise}.
\end{array}  \right.
\]

We now recall that a \emph{system of relations} $R$ for $C=kQ_C/I$ is
a subset $R$ of $\cup_{\ell,h=1}^n e_\ell I e_h$, such that $R$, but no
proper subset of $R$, generates $I$ as an ideal of $kQ_C$.

\begin{lem}\label{lem 1.4}
Let $C$ be a tilted algebra and $R$ be a system of relations for
$C=kQ_C/I$. The quiver $Q_\Crep$ of the cluster repetitive algebra
$\Crep$ is constructed as follows: 
\begin{itemize}
\item[(a)] $ (Q_{\Crep})_0=\{(\ell,i)\mid 1\le \ell\le n, i\in \mathbb{Z}\}$.
\item[(b)] For $(\ell,i),(h,j) \in (Q_{\tilde C})_0$, the set of arrows from $(\ell,i)$ to $(h,j)$ equals \begin{itemize}
\item[(i)] The set of arrows from $\ell$ to $h$ in $Q_C$ if $i=j$,
\item[(ii)] $\textup{Card}\,(R\cap e_h I e_\ell)$ arrows if $i=j+1$,
\end{itemize}   
and there are no other arrows.
\end{itemize}   
\end{lem}

\begin{pf} This follows at once from the above comments and \cite[Theorem 2.6]{ABS}.
\qed
\end{pf}  

Thus the quiver of $\Crep$ can be thought of as consisting of infinitely many copies $(Q_{C_i})_{i\in\mathbb{Z}}$ of the quiver of $C$, joined together by additional arrows from $Q_{C_{i+1}}$ to $Q_{C_i}$, corresponding to $ \Ext^2_C(DC,C)$. In particular, the quiver $Q_\Crep$ is connected if and only if the tilted algebra $C$ is not hereditary.

\begin{example}\label{ex 1.4}
Let $C$ be given by the quiver
\[
\xymatrix@C=60pt@R=30pt{1\ar[r]^\za&2\ar[r]^\zb&3\ar[r]^\zg&4\ar[r]&5}
\]
bound by the relation $\za\zb\zg=0$. 
Then $\Crep$ is given by the infinite quiver
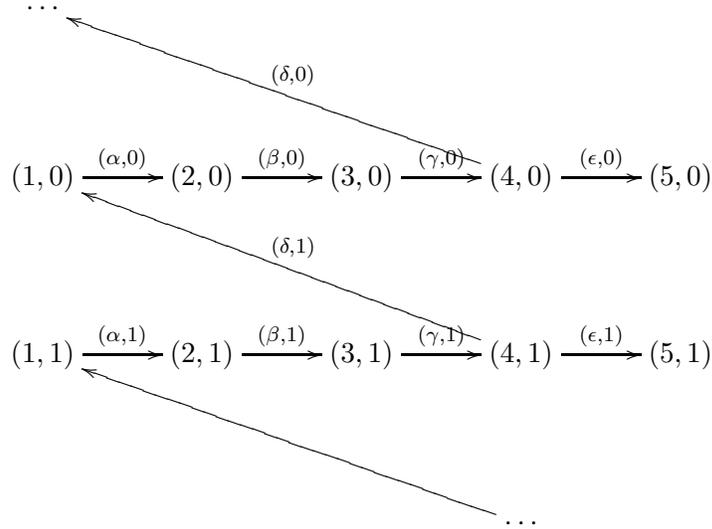
\[
\xymatrix@C=30pt@R=50pt{
\cdots\\
(1,0)\ar[r]^{(\za,0)}&(2,0)\ar[r]^{(\zb,0)}&(3,0)\ar[r]^{(\zg,0)}&(4,0)\ar[r]^{(\ze,0)}\ar[ulll]_{(\zd,0)}&(5,0)\\
(1,1)\ar[r]^{(\za,1)}&(2,1)\ar[r]^{(\zb,1)}&(3,1)\ar[r]^{(\zg,1)}&(4,1)\ar[r]^{(\ze,1)}\ar[ulll]_{(\zd,1)}&(5,1)\\
&&&\cdots\ar[ulll]
}
\]
 bound by the relations 
\[\begin{array}{cc}
(\za,i)(\zb,i)(\zg,i)=0,& (\zd,i+1)(\za,i)(\zb,i)=0,\\
 (\zg,i+1)(\zd,i+1)(\za,i)=0,& (\zb,i)(\zg,i)(\zd,i)=0, \end{array}\] for all $i\in\mathbb{Z}$.
\end{example} 

\end{subsection} 

\end{section}


\begin{section}{The relation with the derived category }\label{sect 2}

\begin{subsection}{}\label{sect 2.1}

Throughout this paper, we let $A$ be a finite dimensional hereditary
algebra, $T$ be a tilting $A$-module and $C=\textup{End}_AT$ be the
corresponding tilted algebra. By Lemma \ref{lem 1.3}, the natural
functor \[\HD\]carries $\DA$ into the category $\textup{Mod}\, \Crep$
of (not necessarily finitely generated) $\Crep$-modules. Since, for
every indecomposable object $X$ in $\DA$, we have $\dim_k\HDX{X}
<\infty$, then its image lies in $\textup{mod}\, \Crep$.

\begin{prop}\label{prop 2.1}
The functor 
\[\HD:\DA\to \textup{mod}\, \Crep\]
 is full and dense and it induces an equivalence of categories 
\[ \DA/\iadd\{\tau F^i T\}_{i\in\mathbb{Z}}  \stackrel{\isomorphe}{\longrightarrow} \textup{mod}\, \Crep ,
\]
where $\iadd\{\tau F^i T\}_{i\in\mathbb{Z}} $ denotes the ideal of $\DA$ consisting of all morphisms which factor through $\add\{\tau F^i T\}_{i\in\mathbb{Z}} $.
\end{prop}

\begin{pf}
We first claim that \[\textup{Ker}\,\HD =\iadd\{\tau F^i
T\}_{i\in\mathbb{Z}}.\] Indeed, let $f:X\to Y$ be a morphism in $\DA$
such that \[\HDX{f}=0.\] By definition of the cluster category, this
means that the induced morphism 
\[\Hom_\CA(\pi T,\pi f):\Hom_\CA(\pi T,\pi X)\to \Hom_\CA(\pi T,\pi Y)
\]
is zero. Therefore $\pi f$ lies in the kernel of $\Hom_\CA(\pi T,-)$, that is, $\pi f$ factors through an object of $\add \pi(\tau T)$. But this implies that $f$ factors through $\add\{\tau F^i T\}_{i\in\mathbb{Z}} $.

Conversely, we prove that any morphism which factors through
$\add\{\tau F^i T\}_{i\in\mathbb{Z}} $ has a zero image. For this, it
suffices to show that the image of any object of the form $\tau F^j T$
(with $j\in\mathbb{Z}$)
is zero. But now 
\[ \HDX{\tau F^j T}=\oplus_{i\in\mathbb{Z}} \Hom_\DA(T,F^{j-i}\tau T) =0
\]
because of Lemma \ref{lem 1.2}.

We now claim that the functor $\HD$ induces an equivalence between
$\add\{ F^i T\}_{i\in\mathbb{Z}}$ and the subcategory $\textup{proj}\, \Crep$ of
$\textup{mod}\, \Crep$ consisting of the projective
$\Crep$-modules. Indeed, by Lemma \ref{lem 1.3}, the
restriction of $\HD$ to $\add\{ F^i T\}_{i\in\mathbb{Z}}$ maps into
$\textup{proj}\, \Crep$. Since, conversely, an indecomposable projective
$\Crep$-module $\check P_0$ is an indecomposable summand of
$\Crep_\Crep=\End(\oplus_{i\in\mathbb{Z}} F^iT)$, then there exists an
indecomposable object $T_0\in\add\{ F^i T\}_{i\in\mathbb{Z}}$ such
that $\check P_0\isomorphe \HDX{T_0}$, that is, the functor is
dense. By Yoneda's lemma, it is full. Let thus $f:T_1\to T_0$ be a
morphism in $\add\{ F^i T\}_{i\in\mathbb{Z}}$ such that
$\HDX{f}=0$. Then $f$ factors through an object of $\add\{\tau F^i
T\}_{i\in\mathbb{Z}}$. Therefore, by Lemma \ref{lem 1.2}, $f=0$. Thus
the functor is faithful and our claim is established. 

It remains to show that the functor $\HD:\DA\to \textup{mod}\, \Crep$ is full and dense. Let $L\in \textup{mod}\, \Crep$ and consider the minimal projective presentation
\[
\xymatrix@C=60pt@R=30pt{
\check P_1\ar[r]^u & \check P_0 \ar[r]& L \ar[r] &0  
}
\]
in $\textup{mod}\, \Crep$. By our claim above, there exist $T_0,T_1\in \add\{F^iT\}_{i\in\mathbb{Z}} $ and a morphism $v:T_1\to T_0$ such that $\HDX{v}=u$. Applying the functor $\HD$ to the triangle
\[
\xymatrix@C=60pt@R=30pt{
T_1\ar[r]^v & T_0 \ar[r] &\overline{L} \ar[r]&T_1[1]  
}
\]
 and using that $T_1[1]\in \add\{\tau F^iT\}_{i\in\mathbb{Z}}$ (because $T_1\in \add\{F^iT\}_{i\in\mathbb{Z}}$ and $F=\tau^{-1}[1]$) yields an exact sequence
\[
\xymatrix@C=50pt@R=30pt{
\check P_1\ar[r]^u & \check P_0 \ar[r]& \HDX{\overline{L}} \ar[r] &0  
}
\]
 in $\textup{mod}\, \tilde C$. Therefore, $L\isomorphe \HDX{\overline{L}}$ and the functor is dense.

Finally, we show that it is full. Let $f:L\to M$ be a morphism in $\textup{mod}\, \tilde C$. Taking minimal projective presentations of $L$ and $M$, we deduce a commutative diagram with exact rows
\[
\xymatrix@C=60pt@R=30pt{
\check{P}_1\ar[d]^{f_1}\ar[r]^u & \check{P}_0\ar[d]^{f_0} \ar[r]& L\ar[d]^{f} \ar[r] &0  \\
\check{P}_1'\ar[r]^{u'} & \check{P}_0' \ar[r]& M \ar[r] &0  \\
}
\]
in $\textup{mod}\, \Crep$. Considering the morphisms $v:T_1\to T_0$ and $v':T_1'\to T_0'$ in 
$\add\{F^i T\}_{i\in\mathbb{Z}} $ corresponding to $u,u'$, respectively, we find a diagram in $\DA$ where the rows are triangles

\[
\xymatrix@C=60pt@R=30pt{
T_1\ar[d]^{g_1}\ar[r]^v & T_0\ar[d]^{g_0} \ar[r]&
\overline{L}\ar@{.>}[d]^{g} \ar[r] &T_1[1]\ar[d]^{g_1[1]}  \\ 
T_1'\ar[r]^{v'} & T_0' \ar[r] &\overline{M} \ar[r] &T_1'[1] ,
}
\]
that is, there exists $g:\overline{L}\to \overline{M}$ such that the
above diagram commutes. Consequently, $\HDX{g}=f$ and the proof is
complete. 
\qed
\end{pf}  
\end{subsection}

\begin{subsection}{}\label{sect 2.2}
It is well-known (see \cite{KR}) that the cluster-tilted algebra $\tilde C$ is 1-Gorenstein, that is, such that for every injective $\tilde C$-module $I$, we have $\pd I\le 1$ and, for every projective $\tilde C$-module $P$, we have $\id P\le 1$. This property clearly lifts to its Galois covering $\Crep$. This also follows from Proposition \ref{prop 2.1}.

\begin{cor}\label{cor 2.2}
The cluster repetitive algebra $\Crep$ is 1-Gorenstein. In particular, $\textup{gl.dim.}\tilde C\in\{1,\infty\}$.
\end{cor}  

\begin{pf}  By \cite[(IV.2.7) p.115]{ASS}, we need to prove that
  \[\Hom_{\Crep}(D\Crep,\tau_{\Crep}I)=0,\]
 for every injective $\Crep$-module $I$. Now under the equivalence of Proposition \ref{prop 2.1}, every injective $\Crep$-module is the image of an object of the form $\tau^2 T_0\in
  \DA$, where $ T_0\in \add\{F^iT\}_{i\in\mathbb{Z}} $.
It thus suffices to show that
\[ \Hom_{\DA}(\oplus_{i\in\mathbb{Z}} \,\tau^2 F^i T,\tau^3 T_0)=0.\]
 But this follows from the fact that
  $\tau$ is an equivalence in $\DA$ and from Lemma \ref{lem 1.2}. Thus, $\Crep$ is 1-Gorenstein. The proof of the second statement is standard (see, for instance, \cite{KR}).
\qed
\end{pf}  
\end{subsection}

\begin{subsection}{}\label{sect 2.3}
The following Lemma is a ``derived'' version of the projectivisation procedure of \cite[II.2.1]{ARS}.

\begin{lem}\label{lem 2.3}
Let $T_0\in\add\{F^iT\}_{i\in\mathbb{Z}} $ and $X\in\DA$, then the map $f\mapsto \HDX{f}$ induces an isomorphism
\[
\Hom_{\DA}(T_0,X)
\isomorphe \Hom_{\Crep}(\Hom(\oplus F^iT,T_0),\Hom(\oplus F^iT,X)).\]
\end{lem}  
\ \\
\begin{pf} Since the surjectivity follows from the fact that the functor $\HD$ is full (see Proposition \ref{prop 2.1}), we prove the injectivity. Assume $\HDX{f}=0$, then $f$ factors through an object of $\add\{\tau F^i T\}_{i\in\mathbb{Z}} $. We then infer from Lemma \ref{lem 1.2} that $f=0$.
\qed
\end{pf}  

\end{subsection}

\begin{subsection}{}\label{sect 2.4}
We now prove the main result of this section.

\begin{thm}\label{thm 2.4}
There exists a commutative diagram of dense functors

\[
\xymatrix@C=80pt@R=40pt{
\DA \ar[r]^\HD \ar[d]_\pi & \textup{mod}\, \Crep
\ar[d]_{G_\zl} \\
\CA \ar[r]^{\Hom(\pi T, -)} & \textup{mod}\, \tilde C
}
\]
where $G_\zl$ is the push-down functor associated to the Galois covering $G:\Crep\to\Crep$.
\end{thm}  

\begin{pf}
Since $\pi(\tau F^i T)=\tau_\CA(\pi T)$ for each $i$, we have
\[\pi(\iadd\{\tau F^i T\}_{i\in\mathbb{Z}}) =\iadd \tau_\CA(\pi
T).\]
 Therefore, using Proposition \ref{prop 2.1}, the functor $\pi$
induces a functor $H:\textup{mod}\, \Crep \to \textup{mod}\, \tilde C$
such that the following diagram commutes: 
\[
\xymatrix@C=80pt@R=40pt{
\DA \ar[r]^\HD \ar[d]_\pi & \textup{mod}\, \Crep
\ar[d]_{H} \\
\CA \ar[r]^{\Hom(\pi T, -)} & \textup{mod}\, \tilde C
}
\]
 We must show that $H=G_\zl$. Let $M$ be a $\Crep$-module and set $\tilde M=H(M)$. We must prove that, for every $a\in \tilde C_0$, we have \[\tilde M(a) = \oplus_{x/a} M(x),\]
where the sum is taken over all $x\in \Crep_0$ in the fibre
$G^{-1}(a)$ of $a$. We use the following notation: for $x\in\Crep_0$,
we denote by $\tilde P_x$ (or $\check P_x$) the corresponding
indecomposable projective $\tilde C$-module (or $\Crep$-module,
respectively) and by $\tilde T_x$ the corresponding summand of $\pi
T$. 

By Proposition \ref{prop 2.1}, there exists an object $\overline{M}\in\DA$ such that $\HDX{\overline{M}}=M$, thus we have

\begin{eqnarray*}
\tilde M(a) &\isomorphe& \Hom_{\tilde C}(\tilde P_a,\tilde M)\\
&\isomorphe& \Hom_{\CA}(\tilde T_a,\pi\overline{M}),
\end{eqnarray*}  
 because no morphism from $\tilde T_a$ to $\pi \overline{M}$ factors
 through $\add(\tau\pi T)$. Let thus $T_x\in\DA$ be such that $\pi
 T_x=\tilde T_x$. Using Lemma \ref{lem 1.3}, we have

\begin{eqnarray*}
\tilde M(a) 
&\isomorphe& \oplus_{i\in\mathbb{Z}}\  \Hom_{\DA }(F^i T_a,\overline{M}) \\
&\isomorphe &\oplus_{x/a}\  \Hom_{\DA}(T_x,\overline{M}),
\end{eqnarray*}
and by Lemma  \ref{lem 2.3}, this is isomorphic to 
\[\oplus_{x/a}\  \Hom_\Crep(\HDX{T_x},\HDX{\overline{M}}). 
\]
We obtain
\begin{eqnarray*}
\tilde M(a) 
&\isomorphe &\oplus_{x/a}\  \Hom_\Crep(\check P_x,{M}) \\
&\isomorphe &\oplus_{x/a}\   M(x).
\end{eqnarray*}  
This completes the proof that $H=G_\zl$. 
Finally, $G_\zl$ is dense because so is the composition $\Hom(\pi T,-)\circ\pi$.
 
\qed
\end{pf} 

\end{subsection}

\begin{subsection}{}\label{sect 2.5}
We deduce the relations between the Auslander-Reiten  quivers of
$\Crep$ and $\tilde C$. 
\begin{cor}\label{cor 2.5}
\begin{itemize}
\item[(a)] The push-down of an almost split sequence of $\textup{mod}\, \Crep$ is an almost split sequence of $\textup{mod}\, \tilde C$. 
\item[(b)] The push-down functor induces an isomorphism of the
  quotient $\zG(\textup{mod}\, \Crep)/\mathbb{Z}$ of the
  Auslander-Reiten  quiver of $\Crep$ onto the Auslander-Reiten  quiver
  of $\tilde C$. 
\end{itemize}   
\end{cor}

\begin{pf} This follows from \cite[3.6]{G} using the density of the push-down functor.
\qed
\end{pf}  

\end{subsection}

\begin{subsection}{}\label{sect 2.6}
Finally, the following proposition is an analog of \cite[3.2]{BMR1},
and the proof can be easily adapted from there. We include it here for
convenience.

\begin{prop}\label{prop 2.6}
The almost split sequences in $\textup{mod}\, \Crep$ are induced by the almost split triangles in $\DA$.
\end{prop}  

\begin{pf}
By \cite{AR}, the image under $\HD$ of a left (or right) minimal almost split morphism is left (or right, respectively) minimal almost split. Let $u:E\to M$ be a right minimal almost split epimorphism in $\textup{mod}\, \Crep$. Then there exists a right minimal almost split morphism $g:Y\to Z$ in $\DA$ such that $\HDX{g}=u$. We have an almost split triangle
\[
\xymatrix@C=60pt@R=30pt{
\tau Z=X \ar[r]^f&Y\ar[r]^g&Z\ar[r]&X[1].
}
\]
 Applying $\HD$, we get an exact sequence
\[
\xymatrix@C=50pt@R=25pt{
\HDX{X} \ar[r]^{\hspace{50pt}f^*}&E\ar[r]^u&M\ar[r]&0.
}
\]
Since $f$ is minimal left almost split, then so is $f^*$. In
particular, $f^*$ is irreducible. Since $u\ne 0$, $f^*$ is not an
epimorphism,  thus it is a monomorphism. Therefore $f^*$ factors through $\tau M$. That is $\tau M\isomorphe \HDX{X}$, because $f^*$ is irreducible.
\qed
\end{pf}  

\end{subsection} 
\end{section}


\begin{section}{The relation with the repetitive algebra }\label{sect 3}

\begin{subsection}{}\label{sect 3.1}
We recall from \cite{HW} that the \emph{repetitive algebra}  $\Chat$
of a finite dimensional algebra $C$ is the self-injective locally
finite dimensional algebra without identity
\begin{equation}\nonumber\Chat  \quad = \quad \left[
  \begin{array}{cccccccccc}
  \ddots &&&\ 0\ \\
&\ C_{i-1}\ \\
&Q_i&\ C_i\ \\
&&Q_{i+1}&\ C_{i+1}\ \\
&\ 0\ &&&\ddots
  \end{array} \right]
\end{equation}
where matrices have only finitely many non-zero entries, $C_i=C$
and $Q_i= DC$ for all $i\in \mathbb{Z}$, all the remaining
entries are zero, addition is the usual addition of matrices and the
multiplication is induced from that of $C$, the $C$-$C$-bimodule
structure of $DC$ and the zero maps  $DC\otimes_C DC \to 0$.
The identity maps $C_i\to C_{i-1}$, $Q_i\to Q_{i-1}$ induce the
so-called Nakayama automorphism $\nu$ of \Chat. The orbit category
$\Chat/<\nu>$ is then isomorphic to the trivial extension $T(C)=C\ltimes
DC$ of $C$ by its minimal injective cogenerator $DC$.

The repetitive algebra  is closely related to the derived category: if
gl.dim.$C<\infty$, then $\DC$ is equivalent, as a triangulated
category, to the stable module category $\smod \Chat$, see
\cite[II.4.9]{H}. 

Let now, as in section \ref{sect 2}, $A$ be a finite dimensional hereditary algebra, $T$
be a tilting $A$-module and $C=\textup{End}_AT$. We denote by $\zO^i$ the
$i$-th syzygy of a module. Also, we identify the \Chat-modules $C_0$
and $C$.

\begin{lem}\label{lem 3.1}
The functor $\sHom_\Chat(\oplus_{i\in\mathbb{Z}}\, \tau^{-i}\zO^{-i}
C,-)$ maps $\smod \Chat $ into $\textup{mod}\, \Crep$ and induces an
equivalence \[ \smod
\Chat/\iadd\{\tau^{1-i}\zO^{-i}C\}_{i\in\mathbb{Z}} \isomorphe
\textup{mod}\, \Crep\]
where $\iadd\{\tau^{1-i}\zO^{-i}C\}_{i\in\mathbb{Z}} $ denotes the
ideal of $\smod \Chat$ consisting of all morphisms which factor
through an object of $\add\{\tau^{1-i}\zO^{-i}C\}_{i\in\mathbb{Z}} $.
\end{lem}

\begin{pf} By Proposition \ref{prop 2.1}, the functor $\HD$ induces an
  equivalence between
  $\DA/\iadd\{\tau F^{i}T\}_{i\in\mathbb{Z}}   $ and $\textup{mod}\,
  \Crep$. By \cite[III.2.10 and II.4.9]{H}, we have
\[\DA\isomorphe \DC \isomorphe \smod \Chat.
\]
Also, under these equivalences, we have
\[ \Crep = \End_\DA(\oplus_{i\in\mathbb{Z}} F^i T) \isomorphe
\sEnd_\Chat(\oplus_{i\in\mathbb{Z}} \, \tau^{-i}\zO^{-i}C) ,
\]
the image of $\tau F^{i}T$ is $\tau^{1-i}\zO^{-i}C$ (for any
$i\in\mathbb{Z}$) and also the functor $\HD$ becomes
$\sHom_\Chat(\oplus_{i\in\mathbb{Z}}  \tau^{-i}\zO^{-i}C,-)$. This
implies the statement.  
\qed
\end{pf}  
\end{subsection}

\begin{subsection}{}\label{sect 3.2}
 We now wish to introduce a different realisation of the cluster
 category. Let $C$ be a tilted algebra, then there exists an
 automorphism $F_C:\smod \Chat \to \smod \Chat$ defined by $F_C =
 \tau^{-1}\zO^{-1}$. We define $\CC$ to be the orbit category  of
 $\smod \Chat$ under the action of $F_C$, that is, the objects of
 $\CC$ are the orbits $(F_C^iX)_{i\in\mathbb{Z}} $ of the objects $X$
 of $\smod \Chat$, and the morphism set from
 $(F_C^iX)_{i\in\mathbb{Z}} $ to $(F_C^iY)_{i\in\mathbb{Z}} $ is 
$\oplus_{i\in\mathbb{Z}} \,\sHom_\Chat(X,F^i_C Y)$.

We denote by $\spihat:\smod \Chat\to \CC$ the projection functor.

\begin{lem}\label{lem 3.2}
Let $A$ be a finite dimensional hereditary algebra, $T$ be a tilting $A$-module and
$C=\textup{End}_AT$. Then  there exists an equivalence $\eta:\CA\to\CC$ such
that the following diagram commutes.
\[
\xymatrix@C=100pt@R=15pt{
\CA\ar[dr]^{\ \Hom_\CA(\pi T,-)}\ar[dd]^{\isomorphe}_{\eta}\\
& \textup{mod}\, \tilde C\\
\CC \ar[ur]_{\ \Hom_\CC(\spihat C,-)}
}
\]
Furthermore, $\Hom_\CC(\spihat C,-)$ is full and dense and induces an
equivalence of categories $\CC/\add(\tau \spihat C) \isomorphe
\textup{mod}\, \tilde C$.
\end{lem}  

\begin{pf} By \cite[III.2.10 and II.4.9]{H}, we have an  equivalence
  of triangulated categories 
\[\DA\isomorphe \DC \isomorphe \smod \Chat.
\]
Under these equivalences, the automorphism $F$ goes to $F_C$, and the
object $T$ onto the \Chat-module $C$. Therefore there is an
equivalence $\eta:\CA\to\CC$ making the shown  diagram  commute. The
last statement follows from \cite{BMR1}.
\qed
\end{pf}   

\end{subsection}

\begin{subsection}{}\label{sect 3.3}
Let $p:\textup{mod}\, \Chat\to \smod \Chat$ denote the canonical
projection. We define the functor $\phi:\textup{mod}\, \Chat \to
\textup{mod}\, \Crep$ to be the composition 
\[
\xymatrix@C=80pt@R=30pt{
\textup{mod}\, \Chat \ar[r]^{p} &\smod\Chat
\ar[r]^{\sHom_{\Chat}(\oplus \tau^{-i}\zO^{-i} C,-)} &
\textup{mod}\, \Crep.   
}
\]
Also, we denote by $\hat P_x$ the indecomposable projective
$\Chat$-module corresponding to an object $x\in \Chat_0$.

\begin{lem}\label{lem 3.3}
The kernel $\mathcal{J}$ of $\phi$ consists of all morphisms factoring
through an object of $\add\{\hat P_x\oplus \tau^{1-i}\zO^{-i}
C\}_{x\in\Chat_0,i\in\mathbb{Z}} $.
\end{lem}
\begin{pf} Clearly, all such morphisms lie in the kernel of
  $\phi$. Conversely, let $f:X\to Y$ be a morphism in $\textup{mod}\,
  \Chat$ such that $\phi(f)=0$. Then $p(f)$ factors through an object
  of $\add\{ \tau^{1-i}\zO^{-i}                       
C\}_{i\in\mathbb{Z}} $, that is, there exist $Z\in \add\{
  \tau^{1-i}\zO^{-i}C\}_{i\in\mathbb{Z}} $ and morphisms $f_2:X\to Z$,
  $f_1:Z\to Y$ such that $p(f)=p(f_1)p(f_2)$. Thus $f-f_1f_2\in\ker
  p$, that is, $f-f_1f_2$ factors through a projective-injective
  \Chat-module $\hat P$. Thus there exist morphisms $g_2:X\to \hat P$,
  $g_1:\hat P\to Y$ such that $f-f_1f_2=g_1g_2$. Therefore
  $f=\left[f_1 g_1\right]\left[
\begin{array}{c}f_2\\g_2
\end{array}  \right]$ factors through $Z\oplus \hat P$.
\qed
\end{pf} 

\end{subsection}

\begin{subsection}{}\label{sect 3.4}
Let now $\pihat$ denote the composition
\[
\xymatrix@C=60pt@R=30pt{
\textup{mod}\, \Chat \ar[r]^p &\smod \Chat \ar[r]^{\spihat}&\CC. 
}
\]
We prove finally our main theorem.

\begin{thm}\label{thm 3.4}
There is a commutative diagram of dense functors
\[
\xymatrix@C=80pt@R=40pt{
\textup{mod}\, \Chat \ar[r]^\phi \ar[d]_\pihat & \textup{mod}\, \Crep
\ar[d]^{G_\zl} \\
\CC \ar[r]^{\Hom(\pihat C, -)} & \textup{mod}\, \tilde C
}
\]
Moreover, $\phi$ is full and induces an equivalence of categories
$\textup{mod}\, \Chat/\mathcal{J} \isomorphe \textup{mod}\, \Crep$.
\end{thm}   

\begin{pf} The commutativity of the diagram follows from Theorem
  \ref{thm 2.4} and Lemma \ref{lem 3.2}, where we use 
  the fact that $\spihat C =\pihat C$. The functor $\pihat$ is dense,
  since it is the composition of two dense functors and, similarly,
  $\phi $ is full and dense,
  since it is the composition of two full and dense functors. 
Finally, the stated equivalence follows from Lemma \ref{lem 3.3}
\qed
\end{pf}  
\end{subsection} 
\begin{subsection}{}\label{sect 3.5} 
The relation between the Auslander-Reiten quivers of $\Chat$ and 
$\Crep$  follows from the next statement.

\begin{prop}\label{prop 3.4}
The almost split sequences in $\textup{mod}\, \Crep$ are induced from
the almost split triangles in $\textup{mod}\, \Chat$.

\end{prop}  

\begin{pf}
Similar to the proof of Proposition \ref{prop 2.6}.
\qed
\end{pf}

\begin{example} 
Let $C$ be the tilted algebra of example \ref{ex 1.4}. We illustrate
the Auslander-Reiten quivers of $\Chat$ and $\Crep$
 in Figure \ref{fig 1}. In
the Auslander-Reiten quiver of $\Chat$, the
positions of the projective-injective modules are marked by  diamonds
and the positions of the indecomposable summands of
$\oplus_{i\in\mathbb{Z}} \,\tau^{1-i}\zO^{-i}C$ are marked by circles.
As we see, removing the points corresponding to those modules in the
Auslander-Reiten quiver of $\Chat$ yields exactly the Auslander-Reiten quiver of $\Crep$.
\begin{figure}
\begin{center}
 $$\xymatrix@C=6pt@R=6pt
{
&&&&{\scriptstyle \diamond}\ar[ddr]
\\
&&{\scriptstyle \diamond}\ar[dr]&&&&&&&&&&&&&&
{\scriptstyle \diamond}
\\
&\cdot\ar[dr]\ar[ur]&&
 {\scriptstyle \circ}\ar[dr]\ar[uur]&&
 \cdot\ar[dr]&&
\cdot\ar[dr]&&
{\scriptstyle \circ}\ar[dr]&&
\cdot\ar[dr]&&
{\scriptstyle \circ}\ar[dr]&&
\cdot\ar[dr]\ar[ur]&&
\\
&&
\cdot\ar[ur]\ar[dr]&&
{\scriptstyle \circ}\ar[ur]\ar[dr]&&
\cdot\ar[ur]\ar[dr]&&
\cdot\ar[ur]\ar[dr]&&
\cdot\ar[ur]\ar[dr]&&
\cdot\ar[ur]\ar[dr]&&
{\scriptstyle \circ}\ar[ur]\ar[dr]&&
\cdot
\\
\cdots&&&
\cdot\ar[r]\ar[dr]\ar[ur]&
\cdot\ar[r]&
\cdot\ar[ur]\ar[r]\ar[dr]&
{\scriptstyle \circ}\ar[r]&
\cdot\ar[ur]\ar[r]\ar[dr]&
\cdot\ar[r]\ar[ddr]&
\cdot\ar[ur]\ar[r]\ar[dr]&
\cdot\ar[r]\ar[dddr]&
\cdot\ar[ur]\ar[r]\ar[dr]&
\cdot\ar[r]&
\cdot\ar[ur]\ar[r]\ar[dr]&
\cdot\ar[r]&
\cdot\ar[ur]\ar[r]\ar[dr]&\cdot
&\cdots
\\
&&&&
\cdot\ar[ur]&&
\cdot\ar[ur]&&
\cdot\ar[ur]&&
\cdot\ar[ur]&&
{\scriptstyle \circ}\ar[ur]&&
\cdot\ar[ur]&&
{\scriptstyle \circ}
\\
&&&&&&&&&{\scriptstyle \diamond}\ar[uur]
&&&&&&&&
\\
&&&&&&&&&&&{\scriptstyle \diamond}\ar[uuur]
}$$
\end{center}

\begin{center}
 $$\xymatrix@C=6pt@R=6pt
 {&\cdot\ar[dr]&&
&& \cdot\ar[dr] & &
\cdot\ar[dr]&& 
&&\cdot\ar[dr]&&
&& \cdot\ar[dr] & &
\\
&&
 \cdot\ar[dr]&&&&
 \cdot\ar[ur]\ar[dr]&& 
 \cdot\ar[dr]&&
 \cdot\ar[ur]\ar[dr]&& 
 \cdot\ar[dr]&&&&
 \cdot&& 
\\
\cdots& &&
 \cdot\ar[r]\ar[dr]&
 \cdot\ar[r]&
 \cdot\ar[dr]\ar[ur]&&
 \cdot\ar[dr]\ar[ur]\ar[r]&
 \cdot\ar[r]&
 \cdot\ar[ur]\ar[dr]\ar[r]&
 \cdot\ar[r]&
 \cdot\ar[ur]\ar[r]&
 \cdot\ar[r]&
 \cdot\ar[r]\ar[dr]&
 \cdot\ar[r]&
 \cdot\ar[ur]\ar[r]&\cdot&\cdots
\\
 &&
 &&
 \cdot\ar[ur]&&
 \cdot\ar[ur]&&
 \cdot\ar[ur]&&
 \cdot\ar[ur]&&
&& \cdot\ar[ur]&&
}$$\end{center}
\caption{Auslander-Reiten quivers  of $\Chat$ and $\Crep$}\label{fig 1}
\end{figure}
\end{example}
\end{subsection} 
\end{section}


\begin{section}{Fundamental domains}\label{sect 4}

\begin{subsection}{}\label{sect 4.1}
Let $C$ be a tilted algebra. We define the \emph{cluster duplicated algebra}
$\Cbar$ of $C$ to be the (finite dimensional) matrix algebra
\[ \Cbar=\left[
\begin{array}{cc}
C_0 & 0 \\ E&C_1
\end{array}  
\right],
\]
where $C_0=C_1=C$ and $E=\Ext^2_C(DC,C)$, endowed with the ordinary
matrix addition, and the multiplication induced from that of $C$ and from the
$C$-$C$-bimodule structure of $ \Ext^2_C(DC,C)$.

Clearly, $\Cbar$ is identified to the quotient algebra of $\Crep$ defined
by the surjection 
\[ \Crep \longrightarrow \left[
\begin{array}{cc}
C_0 & 0 \\ E_1&C_1
\end{array}  
\right],
\]
in the notation of section \ref{sect 1.3}. In particular, the quiver
$Q_\Cbar$ of $\Cbar$ is identified to the full subquiver of $Q_\Crep$
defined by the points 
\[
\{(h,0)\mid h\in (Q_C)_0\}\quad \bigcup \quad \{(h,1)\mid h\in (Q_C)_0\}.
\]
Thus, $Q_\Cbar$ is connected if and only if $C$ is not hereditary.

Since the trivial extension $\tilde C= C\ltimes \Ext^2_C(DC,C) $  is a
subalgebra of \Cbar, the inclusion map $\tilde C\to \Cbar$ defines a
functor $\zeta:\textup{mod}\, \Cbar \to \textup{mod}\, \tilde C$ (by
restriction of scalars).

First, we recall that, denoting be $e_0$ and $e_1$ the matrices
\[ e_0=\left[
\begin{array}{cc}
1&0\\0&0
\end{array}  
\right]
\qquad \textup{and} \qquad e_1=\left[
\begin{array}{cc}
0&0\\0&1
\end{array}  
\right]
\]
then any \Cbar-module can be written in the form $M=(U,V,\mu)$ where
$U=Me_0,$ $V=Me_1$ are $C$-modules, and $\mu:U\otimes_{C}
E\to V$ is the multiplication
map $u\otimes x \mapsto ux$ $(u\in U, x\in E)$.

We then define $\xi:\textup{mod}\, \Cbar\to \textup{mod}\, \tilde C$ as follows. For a \Cbar-module $(U,V,\mu)$,
the $\tilde C$-module $\xi(U,V,\mu)$ has the $C$-module structure of
$U\oplus V$  and the multiplication of $(u,v)\in U\oplus V$ by $x\in
E$ is given by
\[ (u,v) \, x =(0\,,\,\mu(u\otimes x)).
\]
Thus, for $(u,v)\in M$ and $\left[
\begin{array}{cc} c&0\\x&c 
\end{array}  \right] \in \tilde C= C\ltimes E$,
\[ (u,v)\, \left[
\begin{array}{cc} c&0\\x&c 
\end{array}  \right]
=(uc, vc+\mu(u\otimes x)).
\] 
We define in the same way the action of $\xi$ on the morphisms: if
$(g,h): (U,v,\mu)\to (U',V',\mu') $ is a \Cbar-linear map, we put
$\xi(g,h)=g\oplus h : U\oplus V \to U'\oplus V'$ as a $C$-linear map,
the compatibility of this definition with the multiplication by
elements of $E$ follows from the fact that $h\mu=\mu'(g\otimes 1)$.

We now give another description of the functor $\xi$. Let $\xi$ be the
canonical embedding functor of $\textup{mod}\, \Cbar$ into
$\textup{mod}\, \Crep$ (which is obtained by ``extending by zeros''):
it is full, exact, preserves indecomposable modules and their
composition lengths. We have the following easy lemma.

\begin{lem}\label{lem 4.1}
$\xi=G_\zl\circ \zeta .$
\end{lem}  

\begin{pf} This is a straightforward calculation.
\qed
\end{pf}  
\end{subsection}

\begin{subsection}{}\label{sect 4.1 bis}
We have the following remark about the global dimension of $\Cbar$.
\begin{lem}\label{lem 4.1 bis}
\textup{gl.dim. }$\Cbar \le 5$.
\end{lem}

\begin{pf}
This follows from \cite[Corollary 4']{PR}.
\qed
\end{pf}

Easy examples show that this is a strict bound (take for instance $C$
given by the quiver
$\xymatrix{\bullet&\bullet\ar[l]_\zb&\bullet\ar[l]_\za}$ bound by $\za\zb=0$).
\end{subsection}

\begin{subsection}{}\label{sect 4.2}
Before stating the main result of this section, we need the following
notation. Let $\zL$ be any finite dimensional $k$-algebra and $M,N$ be
two indecomposable \zL-modules. A \emph{path} 
 from $M$ to $N$ in $\ind \zL$ is
a sequence of non-zero morphisms
\begin{equation}\nonumber
M=M_0\stackrel{f_1\,}{\longrightarrow} M_1 \stackrel{f_2\,}{\longrightarrow}
  \cdots \stackrel{f_t\,}{\longrightarrow} M_t =N
\end{equation}
with all $M_i$ in $\ind \zL$. 
In this situation we
say that $M$ is a {\em predecessor} of $N$ and write $M\le N$ and that
$N$ is a {\em successor} of $M$. 

If $S_1$ and $S_2$ are two sets of modules, we write $S_1\le S_2$
if every module in $S_2$ has a predecessor in $S_1$, every module in
$S_1$ has a successor in $S_2$, no module in
$S_2$ has a successor in $S_1$  and no module in $S_1$ has a
predecessor in $S_2$. The notation $S_1 < S_2$ stands for 
$S_1\le S_2 $ and $S_1\cap S_2=\emptyset$. 

We define  a \emph{fundamental domain} for the functor $G_\zl$ to be a
full convex subcategory $\zO$ of $\textup{mod}\, \Crep$ such that the
restriction 
\[G_\zl :\zO\longrightarrow \ind\tilde C
\]
is bijective on objects, faithful, preserves irreducible morphisms and
almost split sequences.

Let now $\zS$ be a complete slice in $\textup{mod}\, C$. We denote by
$\zS_i$ the images of $\zS$ in $\textup{mod}\, C_i$ under the
isomorphisms $C_i\isomorphe C$, $i\in\mathbb{Z}$.

\begin{thm}\label{thm 4.2}
Let $\zS$ be a complete slice in $\textup{\textup{mod}\,}\, C$. Then
\[\zO=\{M\in \ind \Cbar\mid \zS_0\le M<\zS_1\}
\]
is a fundamental domain for the functor $G_\zl$.
\end{thm}   

\begin{pf} Without loss of generality, we may assume that $T$ is an
  $A$-module and that $\zS=\Hom_A(T,DA)$. Let 
\[\zO_{\mathcal {D}}=\{X\in \ind  \DA\mid DA \le X < F\,DA\}.
\]
By \cite{BMRRT}, $\zO_{\mathcal{D}}$ is a fundamental domain for the
functor $\pi:\DA\to \CA$.

We first claim that the image $\check\zO$ of $\zO$ under the functor
\[\HD\]
is equal to the full subcategory of $\ind \Crep$ defined by 
\[\check\zO=\{\check M\in \ind \Crep\mid \zS_0\le \check M<\zS_1\}.
\]
We have $\Hom_\DA(F^iT,DA)=0$ unless $i=0$, since $T$ is an
$A$-module. Hence
\begin{eqnarray*} \HDX{DA} &=& \Hom_{\DA}(T,DA)\\
&=& \Hom_A(T,DA)\\
&=& \zS_0
\end{eqnarray*}
Similarly, 
\[ \HDX{FDA}=\zS_1 .\]
 By Proposition \ref{prop 2.1}, this shows our claim.

Now, note that $\check\zO$ is a fundamental domain for the functor
$G_\zl$. 
This indeed follows from Theorem \ref{thm 2.4}, because $\zO_{\mathcal{D}}$
is a fundamental domain for $\pi$ and from Corollary \ref{cor 2.5}.

Finally, we prove that $\check\zO =\zO$.
For this it suffices to prove that $\check \zO \subset \textup{mod}\,
\Cbar$. Slices are sincere, thus every simple $C_0$ (or $C_1$)-module
occurs as a simple composition factor in $\add \zS_0$ (or $\add
\zS_1$, respectively). Let $e$ be the sum of all primitive idempotents
of $\Crep$ corresponding to the simple modules in $C_0$ and $C_1$. We
have just shown that $ e\Crep e = \Cbar $ and 
$\Cbar \subset \textup{Supp}\,\check\zO$, where
$\textup{Supp}\,\check\zO$ is the support of $\check\zO$, that is, the
full subcategory of $\Crep$ generated by all the points $x\in\Crep_0$
such that $Me_x\ne 0$ for some $M\in \check\zO$.

Now we show that $  \Cbar = \textup{Supp}\,\check\zO $. Suppose there is
some $M\in \check \zO$ having a composition factor $S_x$ with $x$ not
in $C_1$ or $C_2$. Assume first that $x$ lies in $C_i$, where $i\ge
2$. Then  there is a nonzero morphism $f:M\to I_x$, where $I_x$ is the
indecomposable injective \Crep-module corresponding to $x$. Since
$I_x$ is a successor of $\zS_2$ and $M$ is a predecessor of $\zS_1$,
lifting this map to the derived category yields a nonzero morphism
from a predecessor of $F\,DA$  to a successor of $F^2 DA$, which is
impossible (we have used the fact that the functor \HD is full, by
Proposition \ref{prop 2.1}).
The proof is entirely similar in case $i\le -1$.

We have  shown that the indecomposable objects in $\zO$ and
$\check \zO$ coincide. Let now $X\to Y$ be an indecomposable morphism
in $\check \zO$. Since $X,Y$ are both \Cbar-modules, then this is an
irreducible morphism in $\textup{mod}\, \Cbar$, hence in $\zO$.
This shows that $\Cbar=\textup{Supp}\,\check\zO$, and the theorem follows. 
\qed
\end{pf}  
\end{subsection}

\end{section}

\end{article}
\end{document}